\documentclass[12pt]{article}
\usepackage{amsmath,amssymb}
\usepackage{enumitem}
\providecommand{\keywords}[1]{\textbf{\textit{Keywords---}} #1}

\usepackage{setspace}
\usepackage{amsmath}

\usepackage{algorithm}
\usepackage{algpseudocode}
\usepackage[table]{xcolor}
\usepackage{graphicx}
\usepackage{epstopdf}
\usepackage{subfigure}
\usepackage{amsmath}

\newtheorem{theorem}{Theorem}[section]
\newtheorem{proposition}[theorem]{Proposition}

\def\dx{\textnormal{dx}}

\def\dt{\textnormal{dt}}

\usepackage{hyperref}

\usepackage{array}
\newcolumntype{P}[1]{>{\centering\arraybackslash}p{#1}}
\title{Extracting sparse high-dimensional dynamics from limited data}
\author{Hayden Schaeffer, Giang Tran, Rachel Ward}

\begin{document}

\maketitle

\begin{abstract}
Extracting governing equations from dynamic data is an essential task in model selection and parameter estimation. The form of the governing equation is rarely known a priori; however, based on the sparsity-of-effect principle one may assume that the number of candidate functions needed to represent the dynamics is very small. In this work, we leverage the sparse structure of the governing equations along with recent results from random sampling theory to develop methods for selecting dynamical systems from under-sampled data. In particular, we detail three sampling strategies that lead to the exact recovery of first-order dynamical systems when we are given fewer samples than unknowns.  The first method makes no assumptions on the behavior of the data, and requires a certain number of random initial samples. The second method utilizes the structure of the governing equation to limit the number of random initializations needed.  The third method leverages chaotic behavior in the data to construct a nearly deterministic sampling strategy. Using results from compressive sensing, we show that the strategies lead to exact recovery, which is stable to the sparse structure of the governing equations and robust to noise in the estimation of the velocity. Computational results validate each of the sampling strategies and highlight potential applications.
 \end{abstract}


\section{Introduction}\label{sec:introduction}
Since the scientific revolution in the seventeenth century, scientists have endeavored to extract increasingly sophisticated physical models from experimental and observational data. This process is usually done manually in the sense that a scientist, with a strong expertise in the field, must examine the data to uncover meaningful physical laws. Parsing through data by hand is often time-consuming, expensive, and infeasible -- especially when the dimension of the system grows. However, with advances in computing and machine learning algorithms, automated discovery of physical models and mathematical equations directly from data is becoming increasingly possible. Some applications of data-based modeling include (but are certainly not limited to) weather predictions, controls for fluid flows \cite{schmid2010dynamic, schmid2012decomposition, proctor2016dynamic}, construction of climate trend models \cite{majda2009normal}, and disease control models \cite{anderson1992infectious,keeling2011modeling}.

Over the last decade, data-driven methods have seen strong growth due to the abundance of data and the development of sophisticated analytical tools. In this work, a computational method for identifying high-dimensional differential equations (which we will refer to as the `model') from under-sampled dynamical data is developed. We will restrict our attention to quadratic models, which contain a large subset of potential physical laws and applications; for example, Lorenz-like systems used in atmospheric science, Fisher's equation and related reaction-diffusion systems, Vlasov-Poisson equations from plasma physics, and fluid dynamics models like the Navier-Stokes equations. The governing equations for such systems are often moreover \emph{sparse} quadratic polynomials, in that the dynamics depend on a small number of variables and second-order interactions.  Indeed, this allows us to recover the governing equations from under-sampled data: we only need to observe a number of samples of the data roughly on the order of the sparsity level of the system. Our method learns the governing models via a sparse optimization problem over a large set of potential candidate functions. It will be shown the our method selects the correct governing dynamical system (under certain conditions) even when the size of the candidate set far exceeds the size of the data set. In this way, these methods are one of the first applications of compressive sensing to model extraction for dynamical systems.

Recently, regression based methods have been developed and applied to model selection and parameter estimation of dynamic data. The authors of \cite{bongard2007, schmidt2009} first proposed the use of regression to select physical laws from synthetic and experimental data. In particular, a symbolic regression method was developed that compared computed derivatives of the data to analytic derivatives of trial functions, while controlling for the total number of trial functions selected (as to avoid overfitting).  In \cite{brunton2016}, a sparsity-promoting method was proposed for extracting dynamical systems by comparing the computed velocity to a large set of potential trial functions. A sequential thresholded least-squares algorithm was used to fit a (redundant) set of trial functions (typically in the form of monomials) to the velocity.  In \cite{schaeffer2017learning, rudy2017data} sparsity-based methods were proposed for learning nonlinear partial differential equations from spatio-temporal data. In \cite{schaeffer2017sms}, a non-convex sparse optimization method was developed to identify the underlying dynamical system from noisy data sets using an integrated trial set. In \cite{mangan2017model}, the authors utilized the Akaike information criterion to rank different sparse solutions using the method from \cite{brunton2016}, allowing for the automated selection of different models when varying the method's free parameter. Using a sparse regression method along with the Takens' embedding theorem, a data-driven method was proposed in \cite{brunton2017chaos} to decompose chaotic systems into intermittently forced linear systems. A sparse convex optimization method was proposed for joint outlier detection and model selection in \cite{tran2016exact}, when the observed data is locally corrupted by high variance noise. Unlike the previous work in the literature, which was mostly empirical, the authors of \cite{tran2016exact} proved that the separation between the outliers and the `clean' data can be exactly recovered. The use of $L^1$ minimization has also appeared in various data-driven and scientific computing applications, for example \cite{schaeffer2013sparse,mackey2014,hou2015sparse, RWWeighted, ABW17, PHD16, tran2015penalty,caflisch2013pdes,ozolinvs2013compressedCPW,brunton2016,bright2013compressive}. Other data-driven methods for learning data structure and approximating dynamics include: the proper orthogonal decomposition 
\cite{berkooz1993proper,holmes2012turbulence}, Koopman representations
\cite{giannakis2015data, williams2014data}, diffusion maps
\cite{coifman2006diffusion,coifman2005geometric,nadler2006diffusion,nadler2006diffusion2}, and dynamic mode decomposition
\cite{schmid2010dynamic, schmid2012decomposition}.

These data-based regression methods often use the $\ell^0$ or $\ell^1$ penalty to promote sparsity in the learned models (\textit{i.e.} to select a small number of active candidates from the large set of potential trial functions).  Soft-thresholding (related to the $\ell^1$ penalty) for sparse recovery and denoising was first proposed in \cite{donoho1995noising}. The $\ell^1$ regularized least squares problem (referred to as the least absolute shrinkage and selection operator or LASSO) was introduced in \cite{tibshirani1996} to reconstruct a sparse vector from linear observations. Conditions under which $\ell^1$ penalized problems admit sparse solutions are detailed in \cite{candes2006, candes2007, donoho2006most} and have led to many applications in imaging and signal processing.  Refined conditions for the setting of function interpolation and approximation were developed in \cite{RV08, R10, CP11}.

In this work, we develop a model selection and parameter estimation method for learning quadratic high-dimensional differential equations from under-sampled data. Using results from compressive sensing and sampling theory, we show that given a certain sampling of the initial data, a convex optimization problem can recover the coefficients and the governing equations exactly even when the data is under-sampled. We detail three sampling strategies, depending on the level of knowledge of the data or the governing equation. In particular, one can decrease the number of samples needed by adding assumptions to the system. It is important to note that the methodology presented in this work can be extended to higher-order polynomial systems and, more generally, governing equations which are sparse with respect to any given bounded orthogonal system. Also, this formulation benefits from well-known numerical methods for solving $\ell^1$ penalized problems. \footnote{Our code is available on our github page: \\ \url{https://github.com/GiangTTran/ExtractingSparseHighDimensionalDynamicsFromLimitedData}.}

\section{Problem Statement}\label{sec:problemstatement}
Consider the dynamical variable $x(t)\in \mathbb{R}^n$ ($n\gg 1$) governed by the equation $\dot{x}=f(x)$ with initial data $x(t_0)=x_0$.  Assume that $f(x)$ is a quadratic vector-valued equation in $x$, which can be written component-wise as:
\begin{equation}
\begin{cases}
\dot{x}_1 &= f_1(x_1,\ldots, x_n)\\
\dot{x}_2 &= f_2(x_1,\ldots, x_n)\\
   \ & \vdots \\
\dot{x}_n&= f_n(x_1,\ldots, x_n).
\end{cases}
\end{equation}
The goal is to learn $f_1, \ldots, f_n$, given $x$ and $\dot{x}$. For a given initial condition $x_0$, assume that we can obtain a sequence of discrete measurements $\{ x(t_1), x(t_2), \ldots, x(t_{m-1})\}$ through either simulations or observations; however, the function $f(x)$ is unknown.  The total number of samples along a given trajectory (including the initialization), denoted as $m$, can be small and thus we will refer to these measurements as a single \textit{burst}. The key here is that by using a small number of bursts corresponding to random initial data, we are able to provide conditions on the recovery of the underlying dynamics. Although we will show the construction and results for quadratic systems, it is important to note that it is possible to generalized to high-order systems.

We construct the under-sampled measurements as follows. Let $k$ be the index for a given burst; that is, if we sample the initial data $x(t_0;k)$ from some random distribution, we can observe the $k^{th}$ burst, \textit{i.e} $\{x(t_0;k), x(t_1;k), x(t_2;k), \ldots, x(t_{m-1};k)\}$. The corresponding velocity along the burst, denoted by $\{\dot{x}(t_0;k), \dot{x}(t_1;k), \dot{x}(t_2;k), \ldots, \dot{x}(t_{m-1};k)\}$, is either observed or computed to some level of accuracy. The collection of trial functions corresponding to the $k^{th}$ burst is denoted as:
{\footnotesize{
\begin{equation}
A^{(k)} =\begin{pmatrix}
1 & x_1(t_0;k) & x_2(t_0;k) & \cdots & x_n(t_0;k) & x^2_1(t_0;k) & x_1(t_0;k)x_2(t_0;k) & \cdots &x^2_n(t_0;k)\\
1 & x_1(t_1;k) & x_2(t_1;k) & \cdots & x_n(t_1;k) & x^2_1(t_1;k) & x_1(t_1;k)x_2(t_1;k) & \cdots &x^2_n(t_1;k)\\
1 & x_1(t_2;k) & x_2(t_2;k) & \cdots & x_n(t_2;k) & x^2_1(t_2;k) & x_1(t_2;k) x_2(t_2;k) & \cdots &x^2_n(t_2;k)\\
& & & \cdots & & &  \\
& & & \cdots & & &  \\
& & & \cdots & & &  \\
1 & x_1(t_{m-1};k) & x_2(t_{m-1};k) & \cdots & x_n(t_{m-1};k) &  x^2_1(t_{m-1};k) & x_1(t_{m-1};k)x_2(t_{m-1};k) & \cdots & x^2_n(t_{m-1};k)\\
\end{pmatrix}
\label{eqn:dictionaryburstmonomial}
\end{equation}
}}    
and the velocity matrix is denoted as:
{\small
\begin{equation}
V^{(k)}=\begin{pmatrix}
\dot{x}_1(t_0;k)& \dot{x}_2(t_0;k) & \cdots &\dot{x}_n(t_0;k)\\
\dot{x}_1(t_1;k)& \dot{x}_2(t_1;k) & \cdots &\dot{x}_n(t_1;k)\\
\dot{x}_1(t_2;k)& \dot{x}_2(t_2;k) & \cdots &\dot{x}_n(t_2;k)\\
& &  \cdots  &  \\
\dot{x}_1(t_{m-1};k)& \dot{x}_2(t_{m-1};k) & \cdots & \dot{x}_n(t_{m-1};k)
\end{pmatrix} .
\end{equation}
}
Let $c_j$ be the vector of coefficients corresponding to the $j^{th}$ governing equation, $1\leq j \leq n$, and define the coefficient matrix:
\begin{equation}
C= \begin{pmatrix}
| & | &  & | \\
c_1 & c_2 & \vdots & c_n \vspace{0.2cm} \\ 
| & | &  & |
\end{pmatrix}.
\end{equation}
Using the $k^{th}$ burst data leads to the subproblem: find $c$ such that $V^{(k)} = A^{(k)} C$. Next, we combine the data over all bursts $k$ from $1,\ldots, K$  by simply concatenating the burst arrays vertically: 
\begin{equation*}
V = \begin{pmatrix}
V^{(1)} \\
V^{(2)} \\
| \\
V^{(K)} 
\end{pmatrix}
\hspace{0.3in} \text{and} \hspace{0.3in}
A = \begin{pmatrix}
A^{(1)} \\
A^{(2)} \\
| \\
A^{(K)} 
\end{pmatrix}
\end{equation*}
and thus the full problem is to find $C$ such that $V = A C$. Let $N:=\dfrac{n^2+3n+2}{2}$ \, be the number of monomials of $n$ variables up to degree two, then the size of matrix $A$ is $mK \times N$ with $mK<N$ therefore the linear inverse problem is ill-posed. Since the coefficient matrix $C$ is sparse, the inversion could be regularized by solving the non-convex optimization problem:
\begin{equation*}
\min_C \ ||C||_{0}  \ \ \  \text{subject to} \ \ \|AC-V\| \leq \sigma, \end{equation*}
where {$||\cdot ||_0$} is the $\ell^0$ penalty, \textit{i.e.} the number of nonzero terms and the norm {$||\cdot ||$} is the maximum of the $\ell^2$ norm of each row. In this way, the inversion is component-wise separable.  The parameter $\sigma>0$ controls for the error between the computed derivative and the true velocity. For any fixed value of $\sigma$, the general $\ell^0$ penalized problem is NP hard \cite{FR11}; and thus we use the convex relaxation known as the $\ell^1$ basis pursuit method:
\begin{equation*}
 \boxed{(\text{M-BP}_\sigma):\quad\min_C \ ||C||_{1}  \ \ \  \text{subject to} \ \ \|AC-V\| \leq \sigma.}
 \end{equation*}
Note that if $V$ and $A$ are given with high accuracy then we can solve the following:
\begin{align*}
\min_C \ ||C||_{1}  \ \ \  \text{subject to} \ \ AC=V.
\end{align*}

The procedure and optimization problem above will be referred to as the monomial basis pursuit M-BP, since the trial set contains monomials up to degree two. We can also repeat the process using the tensorized quadratic Legendre polynomials, which corresponds to changing Equation~\eqref{eqn:dictionaryburstmonomial} to: 
{\footnotesize
\begin{equation}
A_L^{(k)} =\begin{pmatrix}
1 & \sqrt{3}x_1(t_0;k)  & \cdots & \sqrt{3}x_n(t_0;k) & \sqrt{5}\dfrac{3x^2_1(t_0;k)-1}{2} & 3x_1(t_0;k) x_2(t_0;k) & \cdots \\
1 & \sqrt{3}x_1(t_1;k)  & \cdots & \sqrt{3}x_n(t_1;k) & \sqrt{5}\dfrac{3x^2_1(t_1;k)-1}{2} & 3x_1(t_1;k)x_2(t_1;k) & \cdots \\
1 & \sqrt{3}x_1(t_2;k) & \cdots & \sqrt{3}x_n(t_2;k) & \sqrt{5} \dfrac{3x^2_1(t_2;k)-1}{2} & 3x_1(t_2;k)x_2(t_2;k) & \cdots \\
& &  \cdots & & &  \\
1 & \sqrt{3}x_1(t_{m-1};k)  & \cdots & \sqrt{3}x_n(t_{m-1};k)  & \sqrt{5} \dfrac{ 3x^2_1(t_{m-1};k)-1}{2} & 3x_1(t_{m-1};k)x_2(t_{m-1};k) & \cdots \\
\end{pmatrix} 
\label{eqn:dictionaryburstlegendre}
\end{equation}
}          
The resulting inverse problem will be referred to as the Legendre basis pursuit L-BP, and we write it concisely as:
\begin{align*}
\min_{C_L} \ ||C_L||_{1}  \ \ \  \text{subject to} \ \ A_L {C_L}=V.
\end{align*}
 Note that sparsity of the solution of L-BP can be slightly larger than M-BP, since the `pure' quadratic terms now include a constant term.  In particular, if a component of the governing equation is $s$-sparse with respect to $A$, then it is at most $(s+1)$-sparse with respect to $A_L$. A noise-robust extension of this is given by:
\begin{equation*}
 \boxed{(\text{L-BP}_\sigma):\quad\min_{C_L } \ ||C_L||_{1}  \ \ \  \text{subject to} \ \ \|A_LC_L-V\| \leq \sigma.}
 \end{equation*}

\section{Reconstruction Guarantee}
Compressive sensing theory provides reconstruction guarantees for sparse solutions to ill-posed linear inverse problems via the solution to $\ell_1$ optimization problems like L-BP -- let us recall the basics of this theory.  Consider a general system of linear equations, $y = A x$.  If $A$ is underdetermined, or has fewer rows than columns, it is in general impossible to determine $x$ given only $y$ and $A$.   However, if we know a priori that $x \in \mathbb{R}^{N}$ is $s$-sparse, or only has $s \ll {N}$ non-zero entries, the locations of which are unknown, and if the underdetermined matrix $A$ is suitably \emph{incoherent}, it is possible to recover such an $x$ from only ${M} = O(s \log {N})$ measurements as the unique solution to the $\ell_1$-minimization problem $x = \arg \min_z \| z \|_1 \ \text{subject to} \ A z = y$.    Moreover, if there is noise on the measurements $y = Ax + \eta$, then $x$ is well-approximated by the solution $x^{\#} \in \arg \min_z \| z \|_1 \ \text{subject to} \ \| A z - y \|_2 \leq \sigma$ for $\sigma$ appropriately chosen. Thus, by leveraging sparsity, we can potentially overcome the \emph{curse of dimensionality} in the number of measurements we need to take, which often renders high-dimensional inverse problems intractable.  

In our setting, the matrix of coefficients $C$ (equivalently, each of the $n$ columns {$c_1, c_2, \dots, c_n$}) is $s$-sparse, and so we can ask: under what conditions on the measurement  matrix $A$ (or the Legendre-transformed $A_L$) can we apply compressive sensing results and conclude that $C$ is exactly recovered as the solution to {$\text{M-BP}_\sigma$ (or $\text{L-BP}_\sigma$)}?  As it turns out, if each of the $K$ bursts is initialized at a uniformly random point in $[-1,1]^n$, then the $K\times N$ measurement matrix $A_L$ corresponding to only the $K$ initialization measurements, has precisely the incoherence properties that provide optimal compressive sensing results.  Thus, if each of the coefficient vectors $c_k$ is $s$-sparse, then we only need to measure on the order of $K \sim s \log N$ bursts to recover the coefficients exactly.  Each burst need only be measured long enough to get an accurate approximation to the initial velocity. Informally, the size of the burst should be large enough to get a stable approximation to the velocity via a finite difference approximation.

The remarkable fact that sparse vectors can be exactly recovered from vastly underdetermined linear systems of equations cannot be possible for just \emph{any} underdetermined system of equations $Ax = y$ -- for example, the {$M \times N$} matrix  $A$ consisting of the first {$M$} rows of the {$N \times N$} identity matrix maps all sparse vectors whose support does not intersect the first {$M$} coordinates to its null space, rendering it impossible to distinguish them. Indeed, the matrix $A$ must have the \emph{incoherence property} which implies that its null space only intersects the set of sparse vectors trivially, and is sufficiently ``well-conditioned" over the set of sparse vectors to permit stability to noise.   The Legendre-transformed matrix $A_L$ with burst length $m = 0$ and with initializations $x(t_0;1), \dots, x(t_0;K)$ taken as independent and identically distributed as uniform random variables $[-1,1]^n$ satisfies these requirements; in particular, it has two key properties which permit theoretical results on sparse recovery: 
\begin{itemize}
\item Its rows $\omega_1, \omega_2, \dots, \omega_K$ are independent realizations of an isotropic random variable $\omega$; that is, its covariance matrix is the identity matrix $ \mathbb{E}[\omega \omega^*] = {\bf 1}_{N \times N}$,
\item Its rows are \emph{uniformly bounded}: $\max_{1 \leq i \leq N, 1 \leq j \leq K} | A_L(i,j) |^2 \leq 9$.
\end{itemize}
Using these properties, we may apply Theorem 1.2 from \cite{CP11} (see also Theorem 12.22 in \cite{FR11}) to arrive at the reconstruction guarantee.

\begin{theorem}
\label{thm:main}
Consider a dynamical system $\dot{x}=f(x)$ to be recovered from snapshots \\$x(t_0; k), \dots x(t_{m-1}; k)$ and corresponding velocities $\dot{x}(t_0; k), \dots \dot{x}(t_{m-1};k)$, $k=1,2,\dots, K$.   Assume that each component $f(x) = (f_1(x), f_2(x), \dots, f_n(x))$ is a quadratic vector-valued equation in $x,$ and that each of the $f_k$ has at most $s$ out of $N$ polynomial coefficients non-zero.  Suppose that the initialization $x(t_0;k)$
for each of the $k = 1,2,\dots K$ bursts is chosen independently at random from the uniform distribution over $[-1,1]^n$.  
Suppose the number of bursts satisfies
\begin{align}
\label{eq:numburst}
K &\geq 9 c_* s \log(N) \log( \varepsilon^{-1})
\end{align}
where $c_*$ is a universal constant.
Then with probability $1-\varepsilon$, any particular component of the $n$ governing equations in the system $\dot{x}=f(x)$ is recovered exactly via the polynomial coefficients as the unique solution to (L-BP).

\bigskip

\noindent More generally, under the same conditions as above and with the same probability, if the measured gradient terms are only approximate,
 $$\widetilde{\dot{x}}(t_0;k) = \dot{x}(t_0; k) + \tau_{0;k}, \dots , \, {\widetilde{\dot{x}}(t_{m-1};k) = \dot{x}(t_{m-1}; k) + \tau_{m-1;k}},$$ 
such that $\sqrt{\frac{1}{K}\sum\limits_{k=1}^K |\tau_{0;k}|^2} \leq \eta$, then considering $(\text{L-BP}_{\sigma})$ with matrix $A_L$ consisting of only the initial burst data and with $\sigma = \sqrt{K} \eta$, any particular coefficient vector $c_k$ is approximated by a minimizer $c_k^{\#}$ of  $(\text{L-BP}_{\sigma})$ according to
$$
\| c_k^{\#} - c_k \|_2  \leq c^* \sqrt{s} \eta
$$
where $c^*$ is a universal constant.
\end{theorem}
The proof of Theorem  \ref{thm:main} is deferred to the Appendix.  We note that by the union bound (Boole's inequality), if we observe $K' \geq 9 c_* s \log(N) \log( n\varepsilon^{-1})$ bursts, then with probability $1 - \varepsilon$, the error bounds are guaranteed to hold uniformly for each of the $n$ governing equations. Also, the arguments require a uniform bound on the elements of the matrix $A_L$ restricted to the initial data, which can be chosen to be uniform random variables on $[-1,1]^n$. In general, trajectories of polynomial governing systems can become unbound in finite-time, which would make a uniform bound on $A_L$ not possible. However, with the burst approach, one should expect the trajectory to exist for a short time.

  This sparse recovery result can be generalized in several ways; namely,
\begin{itemize}
\item The initializations should be independent and identically distributed random variables, but it is not crucial that they be  uniform random variables on $[-1,1]^n$.  Using results from \cite{RWLegendre}, the initializations could instead be taken to be i.i.d. according to the Chebyshev measure on $[-1,1]^n$, any measure interpolating between the Chebyshev measure and the uniform measure. 
\item The cubic case is an important structure that can arise from special symmetries. Although it is not directly considered here, our approach can be applied to extract cubic governing equations from data. Theorem~\ref{thm:main} holds if we modify Equation~\eqref{eq:numburst} to:
\begin{align*}
K &\geq 27 c_* s \log(N) \log( \varepsilon^{-1})
\end{align*}
 since one can show that the matrix $A_L$ restricted to uniformly random measurements satisfies the bound $\max_{1 \leq i \leq N, 1 \leq j \leq K} | A_L(i,j) |^2 \leq 27$. However, in the cubic case the number of unknowns $N$ can be quite large since $N ={ n +3 \choose 3}$ compared to the quadratic case where  $N ={ n +2 \choose 2}$. Therefore, the number of basis elements $N$ will be extremely large for high-dimensional cubic systems, rendering the $\ell_1$ reconstruction algorithm quite costly without taking into account additional assumptions on the sparsity structure. Constructing a computational efficient approach under reasonable structural conditions will be explored in a future work.
 \item We could have derived reconstruction guarantees for higher-order polynomial systems.  As with the cubic case, the only difference in the theoretical result is that the uniform bound $\max_{1 \leq i \leq N, 1 \leq j \leq K} | A_L(i,j) |^2$ will grow with the maximal polynomial degree; the constant $9$ in Equation \eqref{eq:numburst} will increase accordingly.    
\item We are not even restricted to consider governing equations which are sparse with respect to the polynomial basis; one could consider a different orthonormal basis such as e.g., sines and cosines which are uniformly bounded, and take the initializations of the bursts as i.i.d. random variables according to the orthogonalization measure (or a measure which is ``close to" the orthogonalization measure) of that basis.
 
\end{itemize}

\section{Ergodicity and the Number of Bursts}

Theorem \ref{thm:main} says that any dynamical system described as a system of sparse quadratic ordinary differential equations can be recovered exactly from a small number of randomly initialized bursts.  In particular, the number of bursts need only scale with the sparsity level of the system, and only logarithmically with the ambient dimension.  This result is true for \emph{any} such dynamical system, independent of the behavior of the trajectories along the bursts.  In this sense, it is a ``worst-case"  result.  In many situations, the data {$x(t_0), x(t_1), \dots, x(t_{m-1})$} along a single burst behaves \emph{chaotically}, mimicking the behavior of a random sequence, and in such cases the number of bursts actually required to achieve exact recover should be far smaller -- as our numerical evidence shows, even a single burst often suffices.

While theoretical results concerning the behavior of high-dimensional dynamical systems have remained elusive, recent large-scale simulation studies, such as \cite{chaos1, chaos2}  demonstrate that high-dimensional dynamical systems described by polynomial systems of equations often exhibit chaotic behavior; in fact, such behavior becomes more and more ``probable" as the dimension of the system increases.   If we are in the regime of chaotic behavior, and if we measure snapshots of the system at the time scale of the chaotic dynamics, then for a smooth function $F$, the sequence {$F(x(t_0)), F(x(t_1)), \dots, F(x(t_{m-1}))$} will satisfy a deterministic form of the law of large numbers -- the so-called Birkhoff Ergodic Theorem -- which says that the time average of the sequence will converge to the space average of $F(x)$ with respect to an underlying invariant measure \cite{Birkhoff31}.  For \emph{strongly} ergodic systems, a stronger Central Limit Theorem holds, and the convergence rate can be quantified.  In particular, the rate of convergence is faster if the \emph{correlation} between successive observed values $F(x_k)$ and $F(x_{k+q})$ (for index $q>0$) is smaller, starting with the works \cite{Ruelle68, Ruelle76, Sinai72, Bowen75}.  There is a rich theory on the decay of correlations for certain classes of low-dimensional dynamical systems, beyond the scope of this article. Morally speaking, if we measure a single burst which has the property that a smooth functional applied to the sequence exhibits a fast decay of correlation, then the resulting matrix $A_L$ will have rows which are only weakly correlated, and thus, should almost fit the theoretical requirements for a compressive sensing result of the form of Theorem \ref{thm:main}.  It is worth noting that a more recent paper \cite{RIPless2} generalizes the results from  \cite{CP11} to accommodate matrices with only weakly-correlated rows.  However, even so, the existing compressive sensing theory does not extend to this situation exactly as the sparse signal to be recovered is not \emph{independent} of the measurement matrix.   Exploring a more precise relationship between the level of ergodicity in the system to recover and its effect on reducing the required number of bursts remains an intriguing direction for future study. 

\section{Sampling Strategies and Computational Results}

In this section, we present three sampling strategies that will lead to exact recovery (with high probability). In particular, Strategy 1 uses the construction from Section~\ref{sec:problemstatement}, Strategy 2 is a reduction of Strategy 1 when one can limit the number of unknowns to an $\ell$-sized neighbor (\textit{i.e.} when $f_j$ only depends on  $x_i$, for $i\in [j-\frac{\ell-1}{2}, j+\frac{\ell-1}{2}]$ for odd $\ell$), and lastly Strategy 3 shows that chaotic trajectories greatly reduce the number of random initializations. It is worth noting that Strategy 1 makes no assumptions on the trajectories $x(t)$ or the support of the monomial basis for $f(x)$, Strategy 2 incorporates information on the support of the monomial basis for $f(x)$ but makes no assumptions on the trajectories $x(t)$, and Strategy 3 uses information about the trajectories $x(t)$. In practice, we expect that some combination of these three approaches could be optimal when dealing with various types of data. Note that although we will show examples with quadratic governing equations,  the number of unknowns, $N$, is non-trivial and grows quickly with dimension.

In all cases, we have under-sampled the data and it is easy to check that standard methods, like the least-squares algorithm, do not produce meaningful results. To validate our approaches, we will apply these strategies to the Lorenz 96 system and a quadratic reaction-diffusion equation. The Lorenz 96 system, introduced in \cite{lorenz1996predictability} as an atmosphere model, contains $n>3$ variables $x_1,\cdots, x_n$ and satisfies
\begin{equation}
\dfrac{\dx_k}{\dt} = -x_{k-2}\, x_{k-1} + x_{k-1} \, x_{k+1} - x_k + F,\quad k=1,\cdots, n.
\label{eqn:lorenz96}
\end{equation}
The constant $F$ is independent of $k$ and $x_0=x_n$ and $x_{n+1}=x_1$.  In both the standard monomial basis and in the Legendre monomial basis, the system is 4-sparse (component-wise). Unless otherwise stated, the number of variables is fixed to be $n = 50$ and the constant is set to $F=8$, where chaotic behavior is expected, see { \cite{lorenz1996predictability}}.

The second system we consider is known as the Fisher's equation, a quadratic reaction-diffusion equation, with reaction term $F(x)=x-x^2$. This equation has applications ranging from population dynamics to combustion physics. Its finite difference discretization of the Fisher's equation with $n$ nodes is given by:
\begin{equation}
\dfrac{\dx_k}{\dt} = x_{k+1}-2x_k + x_{k-1} + \gamma (x_k-x_k^2),\quad k=1,\dots, n,
\label{eqn:reactiondiffusion}
\end{equation}
where the coefficients and time-scale have been re-scaled by the grid-spacing and we impose the periodicity condition: $x_0=x_n$ and $x_{n+1}=x_1$. Reaction-diffusion systems have multi-scale phenomena and can produce traveling waves, for more details see, for example, \cite{mckean1975application,van2003front}. In our examples, the solutions remain bounded since the initial data is sampled in the unit box. Note that when transformed to the Legendre basis, each component of Equation~\eqref{eqn:reactiondiffusion} is $5$-sparse.

Throughout this section, we denote $K$ to be the number of initializations, {$m$} is the size of each burst. For each burst, the time-derivative is approximated using central differencing except at the first and last time-step where we use forward and backward differencing (respectively). The dynamical data is generated by solving the ODEs numerically using Runge-Kutta 45 with a finer tolerance than the sampled time-step $dt$. To solve the L-BP problem we use the $spgl1$ algorithm \cite{van2008probing} for all examples here. Other algorithms that could be successfully applied include: {primal-dual \cite{chambolle2011first}, Douglas-Rachford \cite{eckstein1992douglas}, $SpaRSA$ \cite{wright2009sparse}, or the convex optimization package $cvx$ \cite{grant2008cvx}.}

\medskip

\subsection{Strategy 1: $K \sim c\, s\, \log(N)$ Random Initializations} 

A direct consequence of Theorem \ref{thm:main} is that one can recover the governing equation with $K$ random initializations as long as $K\geq c\, s\, \log(N)$. We first determine the number of uniformly random initializations based on the theorem and collect a burst of size $m$ generated from the system starting with each initialization. The data matrix is concatenated vertically over each burst. The parameters used are specified in each experiment. 

 \begin{figure}[b!]
 \centering
  \includegraphics[width = 2.75 in]{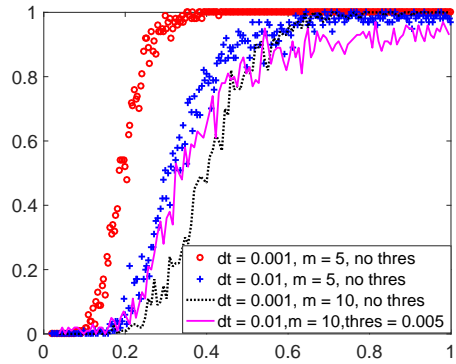}
 \caption{Probability of exact recovery versus the under-sampling rate $K/N$ { for the Lorenz 96 with $n=50$ variables ($N = 1326$), and $F=8$.} For $dt = 0.001$, $K=80$ is needed to achieve $90\%$ probability of success for both $m = 5$ and $m = 10$. For $dt = 0.01$, $K=130$ is needed to achieve $90\%$ probability of success for $m = 5$. }
 \label{fig:Lorenz96_phaseTransition}
 \end{figure}

First, we validate the recovery results by varying the number of initialization $K$ and measuring the probability that the governing system can be recovered. To test this we use the 10th component  of the Lorenz 96 system (see Equation~\eqref{eqn:lorenz96}):
\[\dfrac{dx_{10}}{dt} = -x_{8}\, x_{9} + x_{9} \, x_{11} - x_{10} + 8.\]
 For each $K=5, \ldots, {\dfrac{N}{m}=\dfrac{1326}{m}}$, we repeat the simulation 100 times and record the number of successes to calculate the probability $p$, see Figure \ref{fig:Lorenz96_phaseTransition}. The maximum $K$ is chosen so that we do not oversample the system. From Figure \ref{fig:Lorenz96_phaseTransition}, we can see that for $dt = 0.001$, $K=80$ is needed to achieve $90\%$ probability of success for both {$m = 5$ and $m = 10$}. For a larger time-step, $dt = 0.01$, $K=130$ is needed to achieve $90\%$ probability of success for {$m = 5$}. This is expected since the larger time-steps yield less accurate time-derivative approximations. Note that the recovery degrades when the burst size (in terms of $m$ and $dt)$ increases due to the propagation of error between time-steps. When the data is less correlated in time, the degradation may not occur, see Strategy 3. It is worth noting that although the theoretical universal constant $c_*$ in Equation~\eqref{eq:numburst} has a large upper bound, it is much smaller in practice.

 We also verify the recovery results for Fisher's equation for various $\gamma$. Table~\ref{TableRD} displays the values of the coefficients (as well as their support) when applying the proposed approach to the first component of Equation~\ref{eqn:reactiondiffusion}. The initial data is sampled from the uniform distribution over $[0,1]^n$ and the full data matrix is transformed to $[-1,1]^n$. We set $n=200$ ($N= 20301$), the number of random samples to $159$, and the size of the burst is 5 -- this corresponds to the bound provided in Theorem~\ref{thm:main}. By varying the model parameter $\gamma$, we show that the method selects the correct basis terms and accurately approximates the parameter values (within a few significant digits). It is worth highlighting that as $\gamma$ decreases, the relative scale between the maximum and minimum coefficient increases; however, the recovery results remain relatively stable.

 \begin{table}[t!]
\begin{center}
\rowcolors{1}{lightgray}{lightgray}
 \begin{tabular}{|c||c|c|c|c|c|}\hline
 Terms &   $\gamma=0.25$ & $\gamma=0.1$ &  $\gamma=0.01$&  $\gamma=0$&  true value\\ \hline
   $1$ &$\ 0.000$ &$\ 0.000$  &$\ 0.000$&$\ 0.000$&$\ 0$  \\ \hline
  $x_1$ &$-1.751$ &   $-1.901$ &$\ -1.991$&$\ -1.999$&$\ -2+\gamma$ \\ \hline
  $x_2$ & $1.000$ &   $1.000$ &$\ 1.000$  &$\ 0.999$ &$\ 1$ \\ \hline
   $x_3$ &0 &  0 & 0 & 0 & 0  \\ \hline
          & $\cdots$ &$\cdots$& $\cdots$ & $\cdots$& $\cdots$\\ \hline
 $x_{n-1}$ &0 & 0  &0   & 0  & 0 \\ \hline
$x_n$ & $0.999$ &    $0.999$ &    $1.000$ &    $0.999$&    $1$ \\ \hline
$x_1^2$ & $-0.248$ &   -0.098 &    $-0.008$&    $0$& $-\gamma$\\ \hline
$x_1\,x_2$ & $0$ &    $0$&    $0$ &    $0$&    $0$ \\ \hline
          & $\cdots$ &$\cdots$& $\cdots$ & $\cdots$& $\cdots$ \\ \hline
$x_n^2$ & $0$ &    $0$ &    $0$ &    $0$&    $0$\\ \hline
\end{tabular}
\caption{ Fisher's Equation: Recovery results from under-sampled data generated by Equation~\eqref{eqn:reactiondiffusion}. The results are for the first component of Equation~\eqref{eqn:reactiondiffusion}. As the parameter $\gamma$ varies, the support set of the recovered coefficients remains exact and the corresponding parameters are estimated within a few significant digits. In all cases, we set $n=200$,  $K=159$, and the size of the burst is 5 {($m=5$)}.  }\label{TableRD}
\end{center}
\end{table}

Note that once we have identified the support $S_L$ with respect to the Legendre dictionary matrix, it is possible to ``correct'' the results from Table~\ref{TableRD} by solving:
\begin{align*}
\min_{C}   \|A_M|_{S} C|_{S} - V\|_{2}
\end{align*}
where $A_M|_{S}$ is the monomial dictionary $A_M$ restricted to the corresponding support set $S$.  This debiasing step, when applied to the examples from Table~\ref{TableRD}, produces exact results (up to the fourth significant digit).

\subsection{Strategy 2: $K < c\, s\, \log(N)$ Random Initializations with Localization}

For many ODEs, especially those related to finite dimensional approximations of local PDEs, it is safe to assume that each variable $x_j$ only relates to its $\ell$-neighbors $x_i$, $i\in [j-\frac{\ell-1}{2}, j+\frac{\ell-1}{2}]$ ($\ell$ odd). In particular, the assumption is that the governing equation satisfies:
$$\dot{x}_j = f_j\left(x_{j-\frac{\ell-1}{2}},\ldots, x_ {j+\frac{\ell-1}{2}}\right),$$
for all $j$. With this additional assumption, the number of unknowns in the dictionary matrix is reduced, and thus we can decrease the number of initial conditions needed to guarantee exact recovery. In particular, with the same conditions in Theorem \ref{thm:main}, if additionally it is known that the sparse support is restricted to an $\ell$-sized neighborhood, then if we sample the initial data $K$ times, with {$K\sim c\, s\, \log(\ell)\log(\varepsilon^{-1})$}, then with probability $1-\varepsilon$, the system $\dot{x}=f(x)$ is recovered exactly by the unique solution to (L-BP).

Incorporating this additional condition simply amounts to downsampling the column space of the dictionary {$A$}. In Table~\ref{TableS2}, we consider the $n=1000$ dimensional Fisher's equation, which yields $501501$ unknowns.
By varying the size of the neighborhood, $\ell$, we calculate the minimum number of random initial samples needed for exact recovery. To validate the theoretical scaling, the ratio between the samples needed versus the log of the neighborhood size is shown to be nearly constant. Since the sampling rate is related to the log of the number of unknowns, we see substantial gains between using all terms versus restricting to a neighborhood, but further refinement is not needed. This highlights the benefit of restricting the optimization to an $\ell$-sized neighborhood, even when $\ell$ can only be estimated. One key consequence from this strategy is that the sampling rate is independent of the dimension of the ODE system $n$, allowing this approach to be applied to very large systems.

\begin{table}[t!]
\begin{center}
\rowcolors{1}{lightgray}{lightgray}
 \begin{tabular}{|c||c|c|c|c|c|}\hline
  $\ell$ neighborhood: & All terms & 101 &  51&  31&  11  \\ \hline
  min samples $K$: & 95 & 50 & 43  & 36 & 25   \\ \hline
  ratio $\frac{K}{s\,log(\ell)}$: & $-$ &  2.2 & 2.2  &2.1 & 2.1   \\ \hline
 \end{tabular}
\caption{ Fisher's Equation: For various window sizes $\ell$, the minimum number of random initialization needed for exact recovery is computed for the $n=1000$ dimensional Fisher's equation (number of unknowns is equal to $501501$). The ratio between the samples needed versus the neighborhood size is nearly constant, validating the theory. }\label{TableS2}
\end{center}
\end{table}

\subsection{Strategy 3: Chaotic Systems with $K$ Small and $m$ Large} 

In Strategies 1 and 2, we made no assumption on the behavior of the trajectories $x(t)$.  For Strategy 3, we will assume that the data exhibits chaotic behavior, or uncorrelated long-time behavior. Here we show that if this is the case, we can reduce the number of random initial data needed. In particular, we will show an example of the extreme case where $K=1$, since including additional trajectories with random initial data will only improve the recovery.  Consider collecting data along one trajectory: {
$\{ x(t_0), x(t_1),  \ldots, x(t_{m-1})\}$} with time-steps $dt$ that are  large enough so that $x(t_i)$ and $x(t_{j})$ are sufficiently uncorrelated for $i<j$. The velocity along the trajectory {$\{ \dot{x}(t_0), \dot{x}(t_1), \ldots, \dot{x}(t_{m-1})\}$} is either observed directly (possibly with some error) or calculated by using a fine time-step (smaller than $dt$). In either case, the data is under-sampled, \textit{i.e.}, {$m<N$}. In essence, taking large enough time-steps of a chaotic system mathematically resembles a random ``re-sampling" of the data, thus fitting in with Strategies 1 and 2.

We test this strategy on the first component of the Lorenz 96 system, Equation~\eqref{eqn:lorenz96}, with $n=50$. In Table~\ref{TableOS}, we show that with {500} measurements and $dt = 1.0$, the solution of L-BP identifies the correct terms and approximates the coefficients within the expected error. This shows that it is possible to use the randomness of the data along one trajectory to learn the governing equation. Adding more trajectories with initial data sampled i.i.d. from the uniform distribution, while keeping the total size of the data fixed, only helps the recovery processes.

\begin{table}[h!]
\begin{center}
\rowcolors{1}{lightgray}{lightgray}
 \begin{tabular}{|c||c|c|}\hline
 Terms &   Recovered Coefficients &  True Value\\ \hline
   $ 1$ &$7.986$ &$8$  \\ \hline
  $x_1$ &$-0.9975$ &   $-1.0$ \\ \hline
  $x_2$ & $0.9996$ &   $1.0$ \\ \hline
   $x_3$ &0 &  0   \\ \hline
     $\cdots$     & $\cdots$ &$\cdots$\\ \hline
$x_2\,x_d$ & $0.9996$ &  $1.0$\\ \hline
     $\cdots$     & $\cdots$ &$\cdots$\\ \hline
$x_{d-1}\,x_d$ & $-0.9995$ &    $-1.0$\\ \hline
$x_d^2$ & $0$ &    $0$\\ \hline
\end{tabular}
\end{center}
\caption{Lorenz 96, one trajectory: Using one trajectory of the Lorenz 96 system with $dt = 1.0$, the coefficients are learned using the L-BP. Strategy 3 identifies the correct terms and approximates the coefficients within the expected error.   }
\label{TableOS}
\end{table}

\subsection{Comparisons}
\begin{figure}[h!]
{\includegraphics[width = 2.05 in]{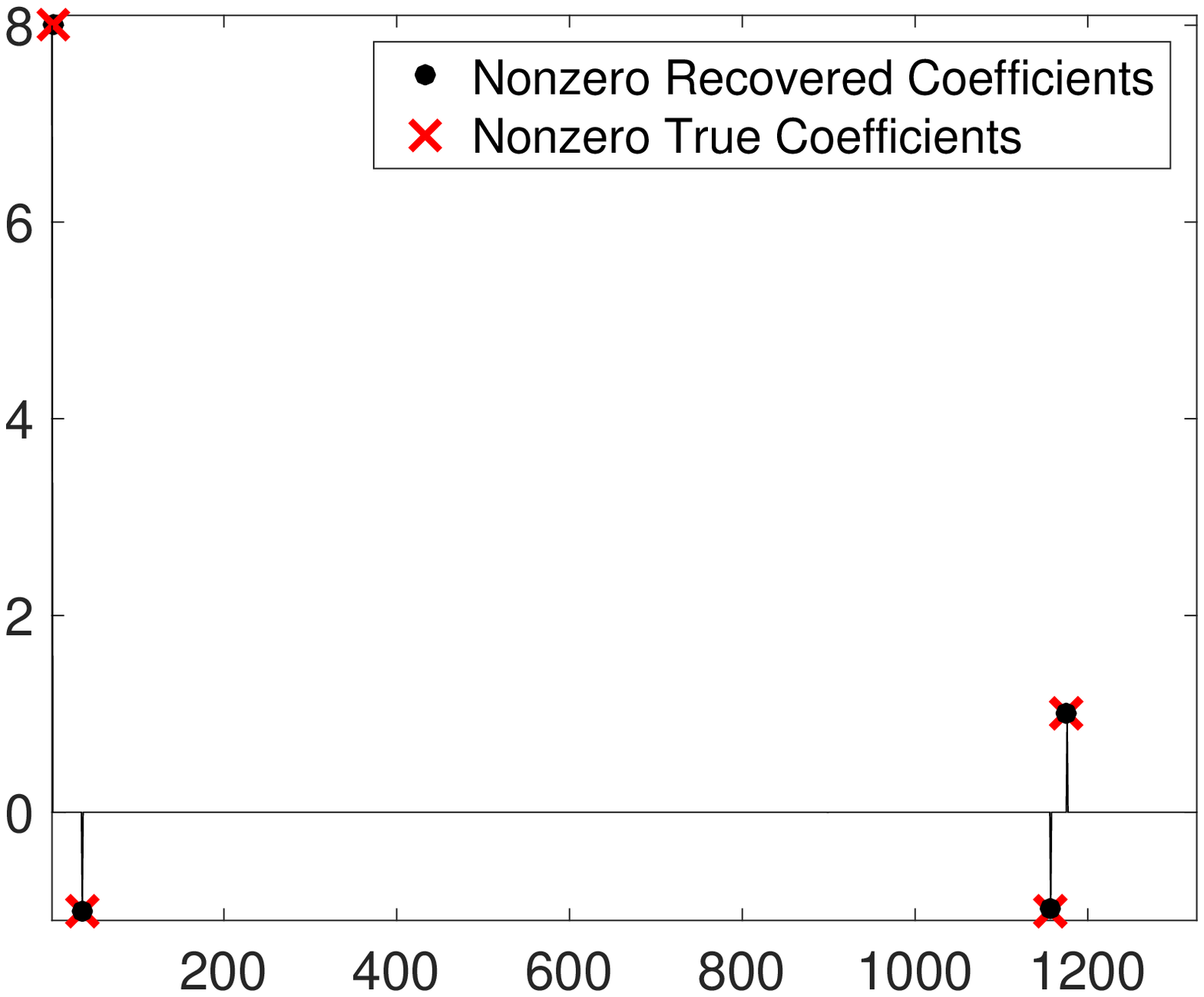}}
{\includegraphics[width = 2.05 in]{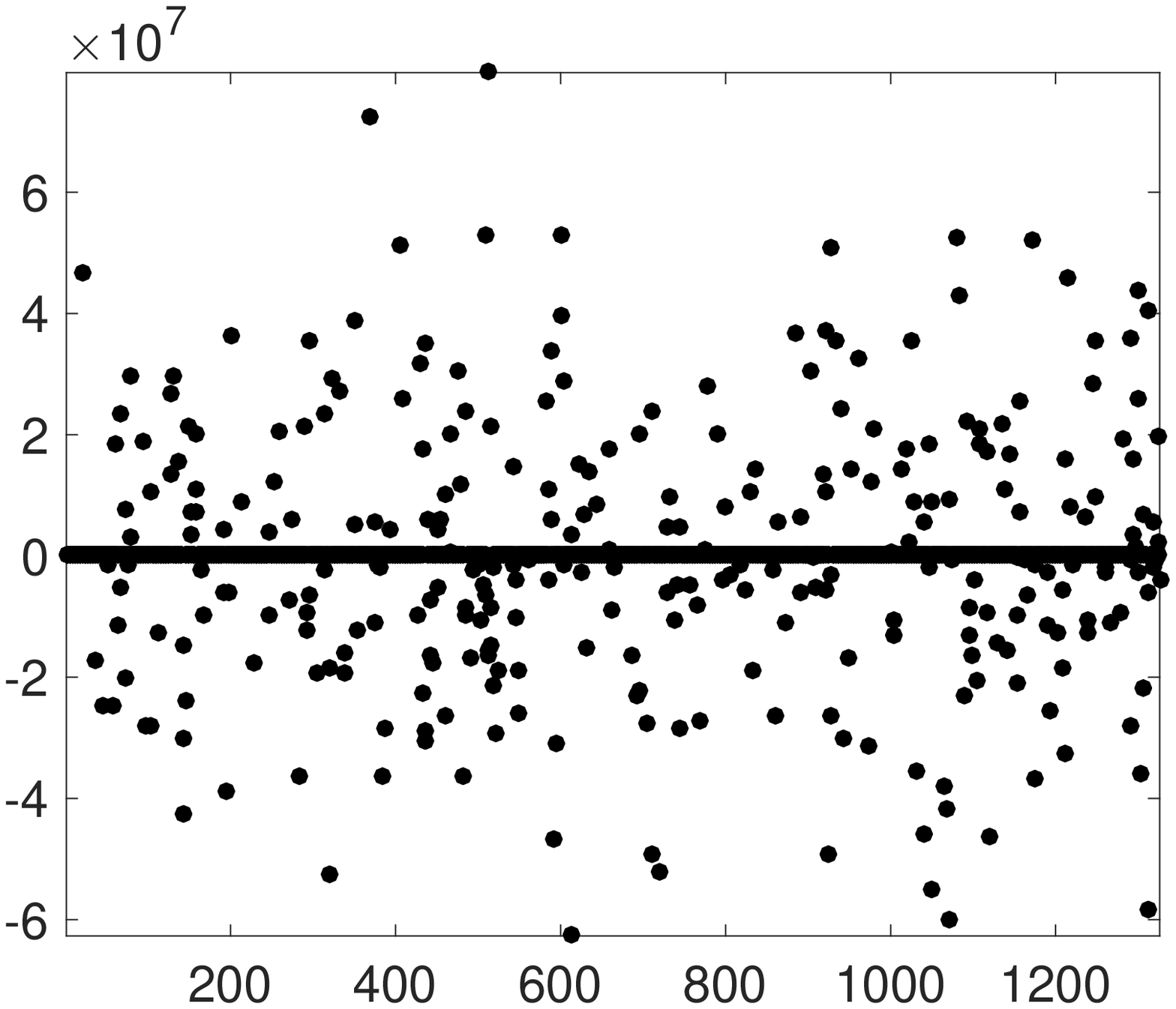}}
{\includegraphics[width = 2.05 in]{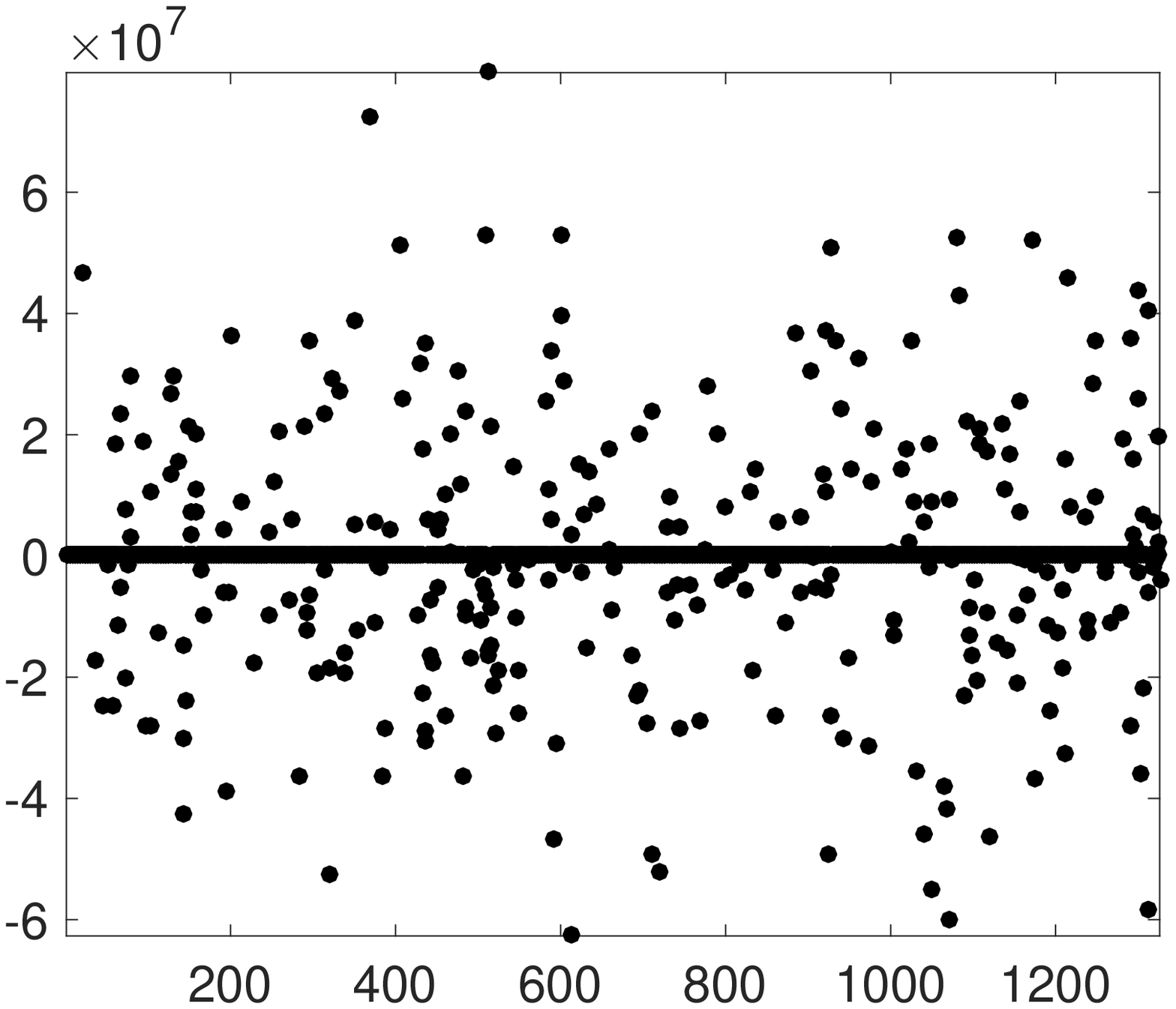}}
\caption{Comparison. The coefficients learned from the L-BP method (left), the least-square algorithm (middle), and  the sequential 
thresholding algorithm (right) for the 35th component of the Lorenz 96 equation with $n=50$, $dt = 0.001$. The threshold parameter for the least-square algorithm and the sequential thresholding algorithm is set to 0.05. }
\label{fig:Lorenz96_compareMethods}
\end{figure}

For comparison, we apply Strategy 1, the standard least-square algorithm, and the linear regression method with sequential thresholding proposed in \cite{brunton2016} on both the Lorenz 96 and Fisher's equation. In Figure \ref{fig:Lorenz96_compareMethods}, the coefficients extracted using the L-BP method (left), the least-square algorithm (middle), and the sequential thresholding algorithm (right) for the 35th component of the Lorenz 96 equation with $n=50$, $dt = 0.001$. The threshold parameter for the least-square algorithm and the sequential thresholding algorithm is set to $\lambda=0.05$. In fact, the solutions of the least-square and the linear regression are the same for this problem. This is the case for all reasonable $\lambda>0$.

Next, we compare the recovery of the coefficients from the first component of the Fisher's Equation~\eqref{eqn:reactiondiffusion}. The model parameter is set to $\gamma=0.1$ and dimension is set to $n=100$, the number of random samples is set to $138$, and the size of the burst is set to 5. The least-squares solution is $653$-sparse and has coefficients on the order of $10^4$. Applying the sequential thresholding algorithm proposed in \cite{brunton2016} with threshold parameter $\lambda \in [5, 5000]$ results in a similar solution. The sparsity matches that of the least-square solution, with $s\in[ 579, 653]$ and has coefficients on the order of $10^4$. Increasing $\lambda$ to achieve a better sparsity level will yield the trivial solution (\textit{i.e.} all zeros). This is likely an effect of their algorithm's dependence on the least-square solution. It is important to note that the algorithm proposed in \cite{brunton2016} was not intended for the case of under-sampling.

\subsection{Measurement Noise}

We consider the effects of noisy state-space measurements on the reconstruction of the sparse coefficient vector. Let $X$ be the data matrix:
\small{
\begin{equation*}
X =\begin{pmatrix}
 x_1(t_0;1) & x_2(t_0;1) & \cdots & x_n(t_0;1) \\
 x_1(t_1;1) & x_2(t_1;1) & \cdots & x_n(t_1;1) \\
& & \cdots &   \\
 x_1(t_{m-1};1) & x_2(t_{m-1};1) & \cdots & x_n(t_{m-1};1) \\
& & \cdots &   \\
 x_1(t_0;k) & x_2(t_0;k) & \cdots & x_n(t_0;k) \\
& & \cdots &   \\
 x_1(t_{m-1};k) & x_2(t_{m-1};k) & \cdots & x_n(t_{m-1};k) \\
 & & \cdots &   \\
   x_1(t_{m-1};K) & x_2(t_{m-1};K) & \cdots & x_n(t_{m-1};K) 
\end{pmatrix}
\end{equation*}
}
and let $Y = X + \eta$  be the matrix of noisy measurements, where $\eta$ is random Gaussian noise. The noise ratio is defined as:
$$\text{Noise Ratio} = \frac{\| X-Y\|_2}{\| X\|_2} \times 100\%$$
and the relative $\ell^2$ error is define as:
$$\text{relative $\ell^2$ error} = \frac{\| c-c_{true}\|_2}{\| c_{true}\|_2} \times 100\%,$$
where $c$ is the computed/learned vector of coefficients and $c_{true}$ is the true coefficients. 
We generate the data from Equation~\eqref{eqn:lorenz96} using $F=8$ and dimension equal to $50$. The parameters are set to $dt=0.001$,  $K=200$, and $m=3$. By adding noise directly to the state-space, we can measure the recovery when the measurements are corrupted. The state variable is corrupted by random Gaussian noise before the derivatives are calculated, thus the value of $V$ and the matrix $A$ will be inaccurate. This makes the problem challenging when the noise is large. The results are summarized in Table~\ref{TableNoise}. As we vary the noise level, we measure the relative $\ell^2$ error and the recovery of the support set.  For noise under $5\%$, the method is stable with respect to the noise. This is consistent with Theorem~\ref{thm:main}. After $5\%$ noise, the $\ell^2$ error jumps and after $6\%$, we cannot reliably recover the support set. In particular, since we know that the true sparsity of Equation~\eqref{eqn:lorenz96} is 4,  we can check if the largest four learned coefficients (in magnitude) coincide with the correct support set. After about $6\%$, the largest four coefficients do not represent the correct support set. This is likely due to the large inaccuracies in $V$, which is not stable to the noise, and inaccuracies in $A$, which scales nonlinearly with the noise.

\begin{table}[b!]
\begin{center}
\rowcolors{1}{lightgray}{lightgray}
 \begin{tabular}{|c||c|c|c|c|}\hline
  Noise Level: & $2.5\%$ &  $5\%$&  $6\%$ &  $7\%$   \\ \hline
  $\ell^2$ error: & $2.6\%$ &   $5\%$  & $10.2\%$ &  $17.2\%$   \\ \hline
  Support Set: &  Y  &  Y & Y &  N  \\ \hline \end{tabular}
\caption{We measure the effects of noise on the recovery of the sparse coefficient vector. For noise under $5\%$, the method is stable to noise. At $5\%$ noise, the $\ell^2$ error jumps and the coefficient vector appears to be polluted by the inaccuracies caused by the noise.}  \label{TableNoise}
\end{center}
\end{table}

\section{Conclusion and Discussion}\label{sec:conclusion}
Extracting dynamical systems remains a difficult task with many open areas of research. Recent work in sparse model selection for dynamical systems has focused on the overdetermined case, where regression must be controlled so as to prevent overfitting. In this work, we utilized the  fundamental idea from compressive sensing to develop several sampling strategies for extracting governing equations from high-dimensional dynamic data. In all cases, the number of measurements is less than the number of unknowns. The main differences between these strategies is the degree of prior knowledge about the data or the governing equations.  If no assumptions on the evolution can be made, then randomizing the initial data is sufficient. If some assumptions on the governing equations are provided, such as locality induced by discretizing a local PDE, then the number of random initial data can be reduced further, to be nearly independent of the dimension of the dataset.  In the third case, if the data is chaotic (or shows a low temporal correlation), then we may reduce the number of initial samples to a fixed number. Using results from compressive sensing, the first two strategies are shown to hold; however, more theory is needed to verify the last case. In several of our experiments,  we have shown that the three strategies are robust to various factors as well as highlighted the benefits of this approach over existing methods. The effective combination of reconstruction guarantees from compressive sensing with sparse learning for dynamical systems presented here opens a wide range of applications. We are currently investigating the use of group sparsity and random sampling approaches for learning dynamic models from multiple data sources. In future work, we also would like to extend the current framework to other bounded orthonormal bases and to learn the dynamics from noisy data.

\section*{Acknowledgments}
The authors would like to thank Amit Singer and Scott McCalla for helpful discussions that improved this manuscript. H.S. acknowledges the support of AFOSR, FA9550-17-1-0125. R.W. and G.T. acknowledge the support of NSF CAREER grant $\#1255631$. 

\section*{Appendix}

For purposes of being as self-contained as possible, we first recall some background on the theory of sparse recovery in \emph{bounded orthonormal systems} via random sampling.  We refer the reader to the text \cite{R10} for more details.

\subsection{Random Sampling in Bounded Orthonormal Systems}

Let ${\cal D} \subset \mathbb{R}^n$ be endowed with a probability measure $\mu$.  Suppose that $\{ \phi_1, \phi_2, \dots, \phi_d \}$ ($d \leq n$) is a (possibly complex-valued) orthonormal system  on ${\cal D}$:
\begin{align}
\int_{{\cal D}} \phi_j(t) \overline{ \phi_k(t)} d\mu(t) = \delta_{j,k} &=  \left\{ \begin{array}{ll} 0 & \text{ if}\ j \neq k \\ 1 & \text{ if}\ j=k \\ \end{array} \right.
\end{align}
We call $\{ \phi_1, \phi_2, \dots, \phi_d \}$ a \emph{bounded orthonormal system} with constant $B \geq 1$ if moreover 
$$
\| \phi_j \|_{\infty} := \text{sup}_{t \in {\cal D}} | \phi_j(t) | \leq B \quad \text{for all}\ j \in 1,2,\dots, d.
$$
Suppose that $t_1, t_2, \dots, t_m \in {\cal D}$ are sampling points which are drawn i.i.d. according to the orthogonalization measure $\mu$, and consider the sampling matrix 
$A \in \mathbb{C}^{m \times n}$ with entries
$$
A_{\ell,k} = \phi_k(t_{\ell}), \quad \ell \in [m], k \in [n].
$$
With high probability, a random matrix formed as such permits stable ``inversion" of the (possibly highly underdetermined) system $y= Ax$ if $x$ is sufficiently sparse, and if ``inversion" is carried out through e.g. solving an $\ell_1$-minimization problem.  The following is a restatement of Theorem 12.22 in \cite{R10}, which is a restatement of a result from \cite{CP11}.

\begin{proposition}
\label{prop:BOS}
Let $x \in \mathbb{C}^n$ and let $A \in \mathbb{C}^{m \times n}$ to be the random sampling matrix associated to a BOS with constant $B \geq 1$.  For $y = Ax + e$ with $\| e \|_2 \leq \eta \sqrt{m}$ for some $\eta \geq 0$, let $x^{\#}$ be a solution to
$$
\min_{z \in \mathbb{C}^n} \| z \|_1 \quad \text{subject to} \ \| Az - y \|_2 \leq \eta \sqrt{m}.
$$
If 
$$
m \geq C B^2 s \log(n) \log(\varepsilon^{-1}),
$$
then with probability at least $1 - \varepsilon$, the reconstruction error satisfies 
$$
\| x - x^{\#} \|_2 \leq C_1 \sigma_s(x)_1 + C_2 \sqrt{s} \eta
$$
where $\sigma_s(x)_1 = \inf_{u: u \text{ is $s$-sparse}} \| x - u \|_1$, and the constants $C, C_1, C_2 > 0$ are universal.
\end{proposition}
We apply this result to prove our reconstruction guarantee.

\subsection{Proof of Theorem \ref{thm:main}}
Recall that $A_{L}$ is the $K \times N$ matrix consisting of those rows corresponding to the $K$ initializations.  The tensor product of univariate Legendre polynomials, normalized as in  the construction of $A_L$, forms a bounded orthonormal system with respect to the uniform measure over $[-1,1]^n$; precisely, $d\mu = \frac{1}{2}dx$.  The Legendre polynomials up to degree 2 are uniformly bounded in magnitude over the domain $[-1,1]^n$ by $3$, as realized by the terms $3 x_j x_k$ at $x_j, x_k = \pm 1$. Thus, $A_{L}$ satisfies the requirements of Proposition \ref{prop:BOS} with $K = 3$.   Moreover, by assumption, each column in the coefficient matrix $C$ associated to the underlying dynamical system consists of at most $s$ nonzero terms.  Thus, after transforming to the Legendre basis, each column in the the Legendre coefficient matrix $C_L$ has at most $s' = s+1$ nonzero terms. Thus, we apply Proposition \ref{prop:BOS} with these parameters to any particular one of the $n$ $\ell_1$-minimization problems. 

The recovery guarantee in the noiseless case $(\eta = 0)$ moreover extends to the optimization problem $\text{L-BP}$ over the \emph{full} set of measurement constraints, not just those constraints corresponding to the burst initializations, since the domain corresponding to the full measurement set is strictly included in the domain corresponding to the reduced measurement set, and the unique minimizer over the reduced measurement set belongs to this subdomain.

\end{document}